\documentclass[11pt]{article}

\usepackage{amssymb,amsmath,amsthm,amsfonts}

\usepackage{epsfig}

\usepackage{rotating}

\usepackage{subfigure}

\newcommand{\pd}[2]{\frac{\partial#1}{\partial#2}}

\begin{document}

\title{Mori-Zwanzig reduced models for uncertainty quantification II: Initial condition uncertainty}
\author{Panos Stinis \\ 
Department of Mathematics \\
University of Minnesota \\
    Minneapolis, MN 55455} 
\date {}

\maketitle

\begin{abstract}
In a recent preprint (arXiv:1211.4285v1) \cite{s12} we addressed the problem of constructing reduced models for time-dependent systems described by differential equations which involve uncertain parameters. In the current work, we focus on the construction of reduced models for systems of differential equations when the initial condition is uncertain. While for both cases the reduced models are constructed through the Mori-Zwanzig formalism, the necessary estimation of the memory parameters is quite different. For the case of uncertain initial conditions we present an algorithm which allows to estimate on the fly the parameters appearing in the reduced model. The first part of the algorithm evolves the full system until the estimation of the parameters for the reduced model has converged. At the time instant that this happens, the algorithm switches to the evolution of only the reduced model with the estimated parameter values from the first part of the algorithm. The viscous Burgers equation with uncertain initial condition is used to illustrate the construction. 
\end{abstract}

\section{Introduction}

The problem of quantifying the uncertainty of the solution of systems of partial or ordinary differential equations has become in recent years a rather active area of research. The realization that more often than not, for problems of practical interest, one is not able to determine the parameters, initial conditions, boundary conditions etc. to within high enough accuracy, has led to a flourishing literature of methods for quantifying the impact that this uncertainty imposes on the solution of the problems under investigation (see e.g. \cite{ghanem,leonenko,ma,nouy,venturi,wan}). However, despite the increase in computational power and the development of various techniques for uncertainty quantification there is still a wealth of problems where reliable uncertainty quantification is beyond reach. The main reason behind the inadequacy is the often high dimensionality (in probability space) of the uncertainty sources. When this uncertainty is coupled with the fact that for practical problems, even solving the corresponding equations for one value of the uncertain parameter (initial condition, boundary condition, $\ldots$) can be very expensive, it results in the uncertainty quantification problem being a rather formidable task. One way to address this problem is to look for reduced models for a subset of the variables needed for a complete description of the uncertainty. 

In the current work, we are concerned with the construction of reduced models for systems of differential equations that arise from polynomial chaos expansions of solutions of a PDE or ODE system. In \cite{s12} we focus on the construction of reduced models when there exists uncertainty in the parameters of the original system. In the current work, we focus on the case that the given PDE or ODE system has uncertain initial condition. Our goal is to construct a reduced model for the evolution of a subset of the polynomial chaos expansions that are needed for a complete description of the uncertainty caused by the uncertainty in the initial condition. There are different methods to construct reduced models for PDE or ODE systems (see e.g. \cite{givon,CS05} and references therein). We choose to use the Mori-Zwanzig (MZ) formalism in order to construct the reduced model \cite{CHK00,CHK3}. 

The main issue with all model reduction approaches is the computation of the memory caused by the process of eliminating variables from the given system (referred to as the full system from this point on) \cite{CS05}. The memory terms are, in general, integral terms which account for the history of the variables that are not resolved. In principle the integrands appearing in the memory terms can be computed through the solution of the orthogonal dynamics equation \cite{CHK3}. However, the solution of this equation is usually very expensive. As we did in \cite{s12}, here too, we utilize a Markovian reformulation of the MZ formalism which allows the calculation of the memory terms through the solution of ordinary differential equations instead of the computation of convolution integrals as they appear in the original formulation. For the example presented in \cite{s12}, the parameters appearing in the reduced model could be estimated {\it without} having to solve the full system. In the current work, this is not possible. However, we present an algorithm which allows the estimation of the necessary parameters on the fly. This means that one starts evolving the full system and use it to estimate the reduced model parameters. Once this is achieved, the simulation continues by evolving {\it only} the reduced model with the necessary parameters set equal to their estimated values from the first part of the algorithm.

Section \ref{mz_formalism} presents a brief introduction to the MZ formalism for the construction of reduced models of systems of ODEs. In Section \ref{memory_comp} we develop the Markovian reformulation of the MZ formalism and show how one can estimate adaptively the parameters appearing in the reduced model. Section \ref{example} applies the reformulation of MZ presented in Section \ref{memory_comp} to the viscous Burgers equation with uncertain initial condition. Finally, in Section \ref{discussion} we discuss certain directions for future work.

\section{Mori-Zwanzig formalism}\label{mz_formalism}

We begin with a brief presentation of the Mori-Zwanzig formalism \cite{CHK00,CHK3}. Suppose we are given the system 
\begin{equation}\label{odes}
\frac{du(t)}{dt} = R (t,u(t)),
\end{equation}
where $u = ( \{u_k\}), \; k \in H \cup G$ 
with initial condition $u(0)=u_0.$ Our goal is to construct a reduced model for the modes in the subset $H.$ The system of ordinary differential equations
we are given can be transformed into a system of  linear
partial differential equations
\begin{equation}
\label{pde}
\pd{\phi_k}{t}=L \phi_k, \qquad \phi_k (u_0,0)=u_{0k}, \, k \in H \cup G
\end{equation}
where $L=\sum_{k \in H \cup G } R_i(u_0) \frac{\partial}{\partial u_{0i}}.$ The solution of \eqref{pde} is
given by $u_k (u_0,t)=\phi_k(u_0,t)$. Using semigroup notation we can rewrite (\ref{pde}) as
$$\pd{}{t} e^{tL} u_{0k}=L e^{tL} u_{0k}$$
Suppose that the vector of initial conditions can be divided as $u_0=(\hat{u}_0,\tilde{u}_0),$ where 
$\hat{u}_0$ is the vector of the resolved variables (those in $H$) and $\tilde{u}_0$ is the vector of the unresolved variables (those in $G$).  Let $P$ be an orthogonal projection on the space of functions of $\hat{u}_0$ and $Q=I-P.$ 

Equation \eqref{pde} 
can be rewritten as 
\begin{equation}
\label{mz}
\frac{\partial}{\partial{t}} e^{tL}u_{0k}=
e^{tL}PLu_{0k}+e^{tQL}QLu_{0k}+
\int_0^t e^{(t-s)L}PLe^{sQL}QLu_{0k}ds, \, k \in H,
\end{equation}
where we have used Dyson's formula
\begin{equation}
\label{dyson1}
e^{tL}=e^{tQL}+\int_0^t e^{(t-s)L}PLe^{sQL}ds.
\end{equation}
Equation (\ref{mz}) is the Mori-Zwanzig identity. 
Note that
this relation is exact and is an alternative way
of writing the original PDE. It is the starting
point of our approximations. Of course, we
have one such equation for each of the resolved
variables $u_k, k \in H$. The first term in (\ref{mz}) is
usually called Markovian since it depends only on the values of the variables
at the current instant, the second is called "noise" and the third "memory". 

If we write
$$e^{tQL}QLu_{0k}=w_k,$$ 
$w_k(u_0,t)$ satisfies the equation
\begin{equation}
\label{ortho}
\begin{cases}
&\frac{\partial}{\partial{t}}w_k(u_0,t)=QLw_k(u_0,t) \\ 
& w_k(u_0,0) = QLx_k=R_k(u_0)-(PR_k)(\hat{u_0}). 
\end{cases} 
\end{equation}
If we project (\ref{ortho}) we get
$$P\frac{\partial}{\partial{t}}w_k(u_0,t)=
PQLw_k(u_0,t)=0,$$
since $PQ=0$. Also for the initial condition
$$Pw_k(u_0,0)=PQLu_{0k}=0$$
by the same argument. Thus, the solution
of (\ref{ortho}) is at all times orthogonal
to the range of $P.$ We call
(\ref{ortho}) the orthogonal dynamics equation. Since the solutions of 
the orthogonal dynamics equation remain orthogonal to the range of $P$, 
we can project the Mori-Zwanzig equation (\ref{mz}) and find
\begin{equation}
\label{mzp}
\frac{\partial}{\partial{t}} Pe^{tL}u_{0k}=
Pe^{tL}PLu_{0k}+
P\int_0^t e^{(t-s)L}PLe^{sQL}QLu_{0k} ds.
\end{equation}

\section{Finite memory}\label{memory_comp}

In this section we describe a reformulation of the problem of computing the memory term which 
does not use the orthogonal dynamics equation. We focus on the case when the memory has a finite extent only. The case of infinite memory is simpler and is a special case of the formulation presented below. Also, the current reformulation allows us to comment on what happens in the case when the memory is very short. 

Let $w_{0k}(t)=P\int_0^t e^{(t-s)L}PLe^{sQL}QLu_{0k} ds=P\int_0^t e^{sL}PLe^{(t-s)QL}QLu_{0k} ds,$ by the change of variables $t'=t-s.$ Note, that $w_{0k}$ depends both on $t$ and the resolved part of the initial conditions $\hat{u}_0.$ We have suppressed the $\hat{u}_0$ dependence for simplicity of notation. If the memory extends only for $t_0$ units in the past (with $t_0 \leq t,$) then $$w_{0k}(t)=P\int_{t-t_0}^t e^{sL}PLe^{(t-s)QL}QLu_{0k} ds.$$ The evolution of 
$w_{0k}$ is given by 
\begin{equation}\label{memory_1}
\frac{dw_{0k}}{dt}=Pe^{tL}PLQLu_{0k}-Pe^{(t-t_0)L}PLe^{t_0 QL}QLu_{0k}+w_{1k}(t),
\end{equation}
where $$w_{1k}(t)=P\int_{t-t_0}^t e^{sL}PLe^{(t-s)QL}QLQLu_{0k} ds.$$ To allow for more flexibility, let us assume that the integrand in the formula for $w_{1k}(t)$ contributes only for $t_1$ units with $t_1 \leq t_0.$ Then $$w_{1k}(t)=P\int_{t-t_1}^t e^{sL}PLe^{(t-s)QL}QLQLu_{0k} ds.$$ 
We can proceed and write an equation for the evolution of $w_{1k}(t)$ which reads
\begin{equation}\label{memory_2}
\frac{dw_{1k}}{dt}=Pe^{tL}PLQLQLu_{0k}-Pe^{(t-t_1)L}PLe^{t_1 QL}QLQLu_{0k}+w_{2k}(t),
\end{equation}
where $$w_{2k}(t)=P\int_{t-t_1}^t e^{sL}PLe^{(t-s)QL}QLQLQLu_{0k} ds.$$ Similarly, if this integral extends only for $t_2$ units in the past with $t_2 \leq t_1,$ then
$$w_{2k}(t)=P\int_{t-t_2}^t e^{sL}PLe^{(t-s)QL}QLQLQLu_{0k} ds.$$
This hierarchy of equations continues indefinitely. Also, we can assume for more flexibility that at every level of the hierarchy we allow the interval of integration for the integral term to extend to fewer or the same units of time than the integral in the previous level. If we keep, say, $n$ terms in this hierarchy, the equation for $w_{(n-1)k}(t)$ will read 
\begin{gather}\label{memory_n}
\frac{dw_{(n-1)k}}{dt}=Pe^{tL}PL(QL)^{n-1}QLu_{0k}- \\
 Pe^{(t-t_{n-1})L}PLe^{t_{n-1} QL}(QL)^{n-1}QLu_{0k}+w_{nk}(t)  \notag
\end{gather}
where $$w_{nk}(t)=P\int_{t-t_n}^t e^{sL}PLe^{(t-s)QL}(QL)^{n}QLu_{0k} ds$$
Note that the last term in \eqref{memory_n} involves the unknown evolution operator for the orthogonal dynamics equation. This situation is the well-known closure problem. We can stop the hierarchy at the $n$th term by assuming that $w_{nk}(t)=0.$

In addition to the closure problem, the unknown evolution operator for the orthogonal dynamics equation appears in the equations for the evolution of $w_{0k}(t),\ldots,w_{(n-1)k}(t)$ through the terms $Pe^{(t-t_0)L}PLe^{t_0 QL}QLu_{0k},\ldots$ $Pe^{(t-t_0)L}PLe^{t_0 QL}(QL)^{n-1}QLu_{0k}$ respectively.

We describe now a way to express these terms involving the unknown orthogonal dynamics operator through known quantities so that we obtain a closed system for the evolution of $w_{0k}(t),\ldots,w_{(n-1)k}(t).$

Since we want to treat the case where $t_0$ is not necessarily small, we divide the interval $[t-t_0,t]$ in $n_0$ subintervals. Define 

\begin{align*}
w_{0k}^{(1)}(t) & =P\int_{t-\Delta t_0}^t e^{sL}PLe^{(t-s)QL}QLu_{0k} ds \\
w_{0k}^{(2)}(t) & =P\int_{t-2 \Delta t_0}^{t- \Delta t_0} e^{sL}PLe^{(t-s)QL}QLu_{0k} ds \\
\ldots & \\
w_{0k}^{(n_0)}(t) & =P\int_{t-t_0}^{t- (n_0-1)\Delta t_0} e^{sL}PLe^{(t-s)QL}QLu_{0k} ds,
\end{align*}
where $n_0 \Delta t_0 = t_0$ and $w_{0k}(t)=\sum_{i=1}^{n_0} w_{0k}^{(i)}(t).$ Similarly, we can define the quantities $w_{1k}^{(1)}(t),\ldots,w_{1k}^{(n_1)}(t)$ 
\begin{align*}
w_{1k}^{(1)}(t) & =P\int_{t-\Delta t_1}^t e^{sL}PLe^{(t-s)QL}QLQLu_{0k} ds \\
w_{1k}^{(2)}(t) & =P\int_{t-2 \Delta t_1}^{t- \Delta t_1} e^{sL}PLe^{(t-s)QL}QLQLu_{0k} ds \\
\ldots & \\
w_{1k}^{(n_1)}(t) & =P\int_{t-t_1}^{t- (n_1-1)\Delta t_1} e^{sL}PLe^{(t-s)QL}QLQLu_{0k} ds,
\end{align*}
where $n_1 \Delta t_1 = t_1$ and $w_{1k}(t)=\sum_{i=1}^{n_1} w_{1k}^{(i)}(t).$ In a similar fashion we can define corresponding quantities for all the memory terms up to  $w_{(n-1)k}(t)=\sum_{i=1}^{n_{n-1}} w_{(n-1)k}^{(i)}(t).$ 

In order to proceed we need to make an approximation for the integrals over the subintervals.

\subsection{Trapezoidal rule approximation}\label{trapezoidal}
We have
\begin{multline*}
w_{0k}^{(1)}(t)  =P\int_{t-\Delta t_0}^t e^{sL}PLe^{(t-s)QL}QLu_{0k} ds   \\
=\biggl[ Pe^{tL}PLQLu_{0k}+Pe^{(t-\Delta t_0)L}PLe^{\Delta t_0 QL}QLu_{0k} \biggr] \frac{\Delta t_0}{2}+ O((\Delta t_0)^3)
\end{multline*}
from which we find
$$Pe^{(t-\Delta t_0)L}PLe^{\Delta t_0 QL}QLu_{0k}=\biggl ( \frac{2}{\Delta t_0} \biggr ) w_{0k}^{(1)}(t) - Pe^{tL}PLQLu_{0k} + O((\Delta t_0)^2)$$
and from \eqref{memory_1}
\begin{equation*}
\frac{dw_{0k}^{(1)}}{dt}=-\biggl ( \frac{2}{\Delta t_0} \biggr ) w_{0k}^{(1)}(t)+ 2Pe^{tL}PLQLu_{0k}+w_{1k}^{(1)}(t)+ O((\Delta t_0)^2).
\end{equation*}
Similarly, for $w_{0k}^{(2)}(t)$ we find
\begin{multline*}
\frac{dw_{0k}^{(2)}}{dt}=\biggl ( \frac{4}{\Delta t_0} \biggr ) w_{0k}^{(1)}(t) \\
-\biggl ( \frac{2}{\Delta t_0} \biggr ) w_{0k}^{(2)}(t) - 2Pe^{tL}PLQLu_{0k} 
+w_{1k}^{(2)}(t)+ O((\Delta t_0)^2)
\end{multline*}
In general,
\begin{multline}\label{memory_1a}
\frac{dw_{0k}^{(i)}}{dt}= -\biggl ( \frac{2}{\Delta t_0} \biggr ) w_{0k}^{(i)}(t) + (-1)^{i+1} 2Pe^{tL}PLQLu_{0k} \\
 +\biggl [  \sum_{j=1}^{i-1}  \biggl ( \frac{4}{\Delta t_0} \biggr ) (-1)^{i+j+1} w_{0k}^{(j)}(t) \biggr ] +w_{1k}^{(i)}(t)+ O((\Delta t_0)^2)  \; \; \text{for}  \; \; i=1,\ldots,n_0.
\end{multline}
Similarly,
\begin{multline*}
\frac{dw_{1k}^{(i)}}{dt}= -\biggl ( \frac{2}{\Delta t_1} \biggr ) w_{1k}^{(i)}(t) + (-1)^{i+1} 2Pe^{tL}PLQLQLu_{0k} \\
 +\biggl [  \sum_{j=1}^{i-1}  \biggl ( \frac{4}{\Delta t_1} \biggr ) (-1)^{i+j+1} w_{1k}^{(j)}(t) \biggr ] +w_{2k}^{(i)}(t)+ O((\Delta t_1)^2)  \; \; \text{for}  \; \; i=1,\ldots,n_1 
\end{multline*} 
$\ldots$
\begin{multline}
\frac{dw_{(n-1)k}^{(i)}}{dt}= -\biggl ( \frac{2}{\Delta t_{n-1}} \biggr ) w_{(n-1)k}^{(i)}(t) + (-1)^{i+1} 2Pe^{tL}PL(QL)^{n-1}QLu_{0k} \\
 +\biggl [  \sum_{j=1}^{i-1}  \biggl ( \frac{4}{\Delta t_{n-1}} \biggr ) (-1)^{i+j+1} w_{(n-1)k}^{(j)}(t) \biggr ] + O((\Delta t_{n-1})^2)  \; \; \text{for}  \; \; i=1,\ldots,n_{n-1}.
\end{multline}
By dropping the $O((\Delta t_0)^2),\ldots, O((\Delta t_{n-1})^2)$ terms we obtain a system of $n_0+n_1+\ldots+n_{n-1}$ differential equations for the evolution of the quantities $w_{0k}^{(1)}(t),\ldots,w_{(n-1)k}^{(n_{n-1})}.$ This system allows us to determine the memory term $w_{0k}(t)=P\int_0^t e^{(t-s)L}PLe^{sQL}QLu_{0k} ds.$ Since the approximation we have used for the integral leads to an error $O(\Delta t)^2,$ the ODE solver should also be $O(\Delta t)^2.$ We have used the modified Euler method to solve numerically the equations for the reduced model. 

Note that the implementation of the above scheme requires the knowledge of the expressions for $Pe^{tL}PLQLu_{0k},\ldots,Pe^{tL}PL(QL)^{n-1}QLu_{0k}.$ Since the computation of these expressions for large $n$ can be rather involved for nonlinear systems (see Section \ref{example}), we expect that the above scheme will be used with a small to moderate value of $n.$ Finally, we mention that the above construction can be carried out for integration rules of higher order e.g. Simpson's rule.

\subsection{Estimation of the memory length}\label{mz_length}

The construction presented above relies on an accurate determination of the memory lengths $\Delta t_0, \Delta t_1,\ldots. \Delta t_{n-1}.$ We present in this section a way to estimate these quantities on the fly. This means that we start evolving the {\it full} system, use it to estimate $\Delta t_0, \Delta t_1,\ldots. \Delta t_{n-1}$ and then switch to the reduced model with the estimated values for $\Delta t_0, \Delta t_1,\ldots. \Delta t_{n-1}.$

For simplicity of presentation we assume that we evolve only $w_{0k}(t)$ and use only one subinterval to discretize the time integrals, i.e. $\Delta t_0=t_0.$ The reduced model reads

\begin{gather}
\frac{d u_{k}}{dt}= Pe^{tL}PLu_{0k} + w_{0k}(t) \label{reduced1} \\
\frac{d w_{0k}}{dt}=2Pe^{tL}PLQLu_{0k} - \frac{2}{t_0} w_{0k}(t)  \label{reduced2}
\end{gather}
for $k \in F. $ 
We can solve \eqref{reduced2} formally and substitute in \eqref{reduced1} to get
\begin{equation}\label{reduced_integral}
\frac{d u_{k}}{dt}= Pe^{tL}PLu_{0k} + \int_0^t e^{-\lambda_0(t-s)}2Pe^{sL}PLQLu_{0k} ds 
\end{equation}
where $\lambda_0=2/t_0.$ Recall that, for the resolved variables, we have from the full system 
\begin{equation}\label{full_split}
\frac{d u_{k}}{dt}= Pe^{tL}PLu_{0k} + Pe^{tL}QLu_{0k}.  
\end{equation}
We would like to estimate the memory decay parameter $t_0$ so that the reduced equation \eqref{reduced_integral} for $u_{k}$ reproduces the behavior of $u_{k}$ as predicted by the full system \eqref{full_split}. We can do that by requiring that the evolution of some integral quantity of the solution is the same when predicted by the reduced and full systems. 

We begin by discretizing the integral term in \eqref{reduced_integral}. Suppose that we are evolving the full system with a step size $\delta t,$ where $t=n_t \delta t$ (note that $n_t$ increases as $t$ increases). If we discretize the integral with the trapezoidal rule we find
\begin{gather}\label{reduced_integral2}
\frac{d u_{k}}{dt}= Pe^{tL}PLu_{0k}  \\
+ [f_{k}(t,\hat{u}_{0})+2 \sum_{j=1}^{n_t-1}e^{-\lambda_0(t-j\delta t)}f_{k}(j\delta t, \hat{u}_{0}) +e^{-\lambda_0t}f_{k}(0,\hat{u}_{0})] \frac{\delta t}{2} \notag
\end{gather} 
where $f_{k}(j\delta t, \hat{u}_{0})=2Pe^{j\delta t L}PLQLu_{0k}$ for $j=0,\ldots,n_t.$ The quantities $f_{k}(j\delta t, \hat{u}_{0})$ can be computed from the full system.  

There is freedom in the choice of the integral quantity whose evolution the reduced model should be able to reproduce. For example, we can use $ \sum_{k \in F}  |u_{k}(t)|^2$ the squared $l_2$ norm of the resolved variables. If we use this integral quantity, then from \eqref{reduced_integral2} and \eqref{full_split} we find that the unknown parameter $t_0$ must satisfy 
\begin{equation}\label{newton1}
 \sum_{k \in F} 2 Re \{  I_{k}(t,t_0) u^*_{k}(t) \} =   \sum_{k \in F}  2 Re \{  Pe^{tL}QLu_{0k} u^*_{k}(t) \} , 
\end{equation} 
where $$I_{k}(t,t_0)= [f_{k}(t,\hat{u}_{0})+2 \sum_{j=1}^{n_t-1}e^{-\lambda_0(t-j\delta t)}f_{k}(j\delta t, \hat{u}_{0}) +e^{-\lambda_0t}f_{k}(0,\hat{u}_{0})] \frac{\delta t}{2}$$ and $Re\{\cdot\}$ denotes the real part.

Let $y=\exp[-\lambda_0\delta t].$ Then, 
\begin{equation}\label{newton2}
I_{k}(t,t_0)= [f_{k}(t,\hat{u}_{0})+2 \sum_{j=1}^{n_t-1} y^{n_t-j}  f_{k}(j\delta t, \hat{u}_{0}) +y^{n_t}f_{k}(0,\hat{u}_{0})] \frac{\delta t}{2}.
\end{equation} 
With this identification, equation \eqref{newton1} becomes a polynomial equation for $y$ with $y \in [0,1].$ It is not difficult to solve equation \eqref{newton1} with an iterative method, for example Newton's method. For the numerical results we present in Section \ref{example}, Newton's method converged to double precision accuracy within 4-5 iterations. After an estimate $\hat{y}$ has been obtained, we can find the estimate $\hat{t}_0$ of $t_0$ (recall $\lambda_0=2/t_0$) from 
\begin{equation}\label{newton3}
\hat{t}_0=-\frac{2 \delta t}{\ln \hat{y}} .
\end{equation}

\subsubsection{Determination of optimal estimate $\hat{t}_0$}\label{mz_optimal}
For each time instant $t$ we can obtain through equations \eqref{newton1} and \eqref{newton3}, an estimate $\hat{t}_0(t)$ for $t_0.$ Thus, the most important issue that we have to address is that of deciding which is the best estimate of $t_0.$ In other words, at what time $t_f$ should we stop estimating the value of $t_0$ so that we can use the estimated value $\hat{t}_0(t_f)$ to evolve the reduced model from then on. 

We define $\epsilon(t)=\underset{l\in [1,n_t]}{\max} |\hat{y}^l (t+\delta t)-\hat{y}^l (t)|.$ The quantity $\epsilon(t)$ monitors the convergence of not only the value of the estimate $\hat{y}$ as a function of the time $t$, but of the whole function $e^{-\lambda_0(t-s)}.$ Ideally, $\epsilon(t)$ converges to zero with increasing $t.$ That will be the case if the approximation of the memory term only through $Pe^{tL}PLQLu_{0kr}$ is enough (see \eqref{reduced1}-\eqref{reduced2}). However, this will not always be the case. If keeping $Pe^{tL}PLQLu_{0kr}$ is not enough, then $\epsilon(t)$ will decrease with increasing $t$ up to some time $t_{min}$ when it will reach a nonzero minimum. After that time, it starts increasing. This signals that keeping only $Pe^{tL}PLQLu_{0kr}$ is {\it not enough} to describe accurately the memory. 

In order to proceed we have two options: (i) construct a higher order model and (ii) identify $t_f=t_{min}$ and thus $\hat{t}_0(t_f)=\hat{t}_0(t_{min}).$ Results for higher order models will be presented elsewhere (see also discussion in Section \ref{discussion}). In the numerical experiments we present in the next section we have chosen $\hat{t}_0(t_f)=\hat{t}_0(t_{min}).$ Note that the procedure just outlined allows the automation of the algorithm. This means that there is no adjustable reduced model parameter that needs to be specified at the onset of the algorithm. 

We are now in a position to state the adaptive Mori-Zwanzig algorithm which constructs a reduced model with the necessary memory term parameter $t_0$ estimated on the fly. 

\vskip14pt
{\bf Adaptive Mori-Zwanzig Algorithm}
\begin{enumerate}
 \item
Evolve the full system and compute, at every step, the estimate $\hat{t}_0(t).$ Use estimates of $t_0$ from successive steps to calculate $\epsilon(t)=\underset{l\in [1,n_t]}{\max} |\hat{y}^l (t+\delta t)-\hat{y}^l (t)|.$  
\item
When $\epsilon(t)$ reaches a minimum (possibly non zero) value at some instant $t_{min}$, pick $\hat{t}_0(t_{min})$ as the final estimate of $t_0.$
\item
For the remaining simulation time, switch from the full system to the reduced model. The reduced model is evolved with the necessary parameter $t_0$ set to its estimated value $\hat{t}_0(t_{min}).$  
\end{enumerate}

This procedure can be extended to the computation of optimal estimates for $t_1,t_2,\ldots,$ i.e. when we evolve, in addition to $w_{0k}(t),$ the quantities $w_{1k}(t),w_{2k}(t),\ldots.$ Results for such higher order models will be presented elsewhere.

\section{Burgers equation with uncertain initial condition}\label{example}

In this section we show how the above MZ formulation can be used for uncertainty quantification for the one-dimensional Burgers equation with uncertain initial condition. The equation is given by 
\begin{equation}\label{burgersequation}
u_t+u u_x = \nu u_{xx},
\end{equation}
where $\nu > 0.$ Equation (\ref{burgersequation}) should be supplemented with an initial condition $u(x,0)=u_0(x)$ and boundary conditions. We solve (\ref{burgersequation}) in the interval $[0,2\pi]$ with periodic boundary conditions. This allows us to expand the solution in Fourier series
$$u_{N}(x,t )=\underset{k \in F}{\sum} u_k(t) e^{ikx},$$
where $F=[-\frac{N}{2},\frac{N}{2}-1].$ The equation of motion for the Fourier mode $u_k$ becomes
\begin{equation}
\label{burgersode}
 \frac{d u_k}{dt}=- \frac{ik}{2} \underset{p, q \in F}{\underset{p+q=k  }{ \sum}} u_{p} u_{q}  -\nu k^2 u_k.
\end{equation}
We assume that the initial condition $u_0(x)$ is uncertain (random) and can be expanded as $u_0(x,\xi)= (\alpha_0+ \alpha_1 \xi) v_0(x)$ where $\xi$ is uniformly distributed in $[-1,1]$ and $v_0(x)$ a given function. In the numerical experiments we have taken $\alpha_0 =\alpha_1=1$ and $v_0(x)=\sin x.$ 

To proceed we expand the solution $u_k(t,\xi)$ for $k \in F$ in a polynomial chaos expansion using Legendre polynomials which are orthogonal in the interval $[-1,1].$ In particular, we have that $$\int_{-1}^1L_i (\xi) L_j(\xi) d \xi=\frac{2}{2i+1} \delta_{ij},$$ where $L_i(\xi)$ is the Legendre polynomial of order $i.$ For each wavenumber $k$ we expand the solution $u_k(t,\xi)$ of \eqref{burgersodemz} in Legendre polynomials and keep the first $M$ polynomials  
\begin{equation}\label{ode_expansion}
u_k(t,\xi)\approx \sum_{i=0}^{M-1} u_{ki}(t) L_i(\xi), \; \; \text{where} \; \; \xi \sim U[-1,1].
\end{equation}
Similarly, the initial condition can be written as $u_0(x,\xi) = \sin x \sum_{i=0}^1\alpha_i L_i(\xi)$ since $L_0 (\xi)=1$ and $L_1 (\xi)=\xi.$
Substitution of \eqref{ode_expansion} in \eqref{burgersodemz}, use of the expansion of the viscosity coefficient and of the orthogonality properties of the Legendre polynomials gives 
\begin{equation}\label{burgersodemz_system}
\frac{du_{kr}(t)}{dt}=- \frac{ik}{2} \sum_{l=0}^{M-1} \sum_{m=0}^{M-1} \underset{p, q \in F}{\underset{p+q=k  }{ \sum}} u_{pl} u_{qm}  c_{lmr} - k^2 u_{kr}
\end{equation}
for $k \in F$ and $r=0,\ldots,M-1.$ Also $$c_{lmr}=\frac{E[ L_l (\xi) L_m(\xi) L_r(\xi) ] }{E[L^2_r(\xi) ]},$$
where the expectation $E[\cdot]$ is taken with respect to the uniform density on $[-1,1].$ The expectation on the denominator of the expression for $c_{lmr}$ is $E[L^2_r(\xi) ]=\int_{-1}^1 L^2_r(\xi) \frac{1}{2}d\xi=\frac{1}{2r+1},$ while the expectation on the numerator can be computed accurately using Gaussian quadrature with Legendre nodes. The Legendre polynomial triple product integral defines a tensor which has the following sparsity pattern: $E[ L_l (\xi) L_m(\xi) L_r(\xi) ]=0,$ if $ l+m < r$ or $l+r < m$ or $m+r < l$ or $l+m+r= \text{odd}$ \cite{gupta}. Due to this sparsity pattern, for a given value of $M$ only about $1/4$ of the $M^3$ tensor entries are different from zero. 

\subsection{MZ reduced model}\label{mz_ode_example}

To conform with 
the Mori-Zwanzig formalism we set 
$$R_{kr}(u)=- \frac{ik}{2} \sum_{l=0}^{M-1} \sum_{m=0}^{M-1} \underset{p, q \in F}{\underset{p+q=k  }{ \sum}} u_{pl} u_{qm}  c_{lmr} -k^2  u_{kr} ,$$
where $u=\{u_{kr}\}$ for $k \in F$ and $r=0,\ldots,M-1.$ Thus, we have
\begin{equation}
\label{burgersodemz}
\frac{d u_{kr}}{dt}=R_{kr}(u) 
\end{equation}
for $k \in F$ and $r=0,\ldots,M-1.$  
We proceed by dividing the variables in resolved and unresolved. In particular, we consider as resolved the variables $\hat{u}=\{u_{kr}\}$ for $k \in F$ and $r=0,\ldots,\Lambda-1,$ where $\Lambda < M.$ Similarly, the unresolved variables are $\tilde{u}=\{u_{kr}\}$ for $k \in F$ and $r=\Lambda,\ldots,M-1.$ In the notation of Section \ref{mz_formalism} we have $H= F \cup (0,\ldots,\Lambda-1)$ and $G= F\cup (\Lambda,\ldots,M-1).$ In other words, we resolve, for all the Fourier modes, only the first $\Lambda$ of the Legendre expansion coefficients and we shall construct a reduced model for them. 

The system (\ref{burgersodemz}) is supplemented by the initial 
condition $u_0=(\hat{u}_0,\tilde{u}_0).$ We focus on initial conditions where 
the unresolved Fourier modes are set to zero, i.e. $u_0=(\hat{u}_0,0).$ We also define $L$ by 
$$L=\sum_{k \in F}\sum_{r=0}^{M-1} R_{kr}(u_0) \frac{\partial}{\partial u_{0kr}}.$$ 
To construct a MZ reduced model we need to define a projection operator $P.$ For a function $h(u_0)$ of all the 
variables, the projection operator we will use is defined by $P(h(u))=P(h(\hat{u}_0,\tilde{u}_0))=h(\hat{u}_0,0),$ i.e. 
it replaces the value of the unresolved variables $\tilde{u}_0$ in any function $h(u_0)$ by zero. Note that this choice of projection is consistent with the initial conditions we have chosen. Also, we define the Markovian term 
$$ PLu_{0k}=PR_k(u_0)=- \frac{ik}{2} \sum_{l=0}^{\Lambda-1} \sum_{m=0}^{\Lambda-1} \underset{p, q \in F}{\underset{p+q=k  }{ \sum}} u_{0pl} u_{0qm}  c_{lmr} -k^2  u_{0kr}.$$ 
The Markovian term has the same functional form as the RHS of the full system but is restricted to a sum over only the first $\Lambda$ Legendre expansion coefficients  for each Fourier mode. 

For the the term $PLQLu_{0kr}$ we find

\begin{equation}\label{burgersmemory1}
PLQLu_{0kr}=2\times \biggl [   - \frac{ik}{2}   \sum_{l=\Lambda}^{M-1} \sum_{m=0}^{\Lambda-1} \underset{p, q \in F}{\underset{p+q=k  }{ \sum}} PLu_{0pl} u_{0qm}  c_{lmr} \biggr ] .
\end{equation}

Finally, to implement any method to solve equation \eqref{newton1} for the estimation of $t_0$ we need to specify the RHS of the equation \eqref{newton1}. This requires the evaluation of the expression $Pe^{tL}QLu_{0kr}.$ For the case of the viscous Burgers equation, we find
\begin{gather}\label{newton4}
Pe^{tL}QLu_{0kr}=2(- \frac{ik}{2}) \sum_{l=\Lambda}^{M-1} \sum_{m=0}^{\Lambda-1} \underset{p, q \in F}{\underset{p+q=k  }{ \sum}} u_{pl} u_{qm}  c_{lmr} \\
- \frac{ik}{2} \sum_{l=\Lambda}^{M-1} \sum_{m=\Lambda}^{M-1} \underset{p, q \in F}{\underset{p+q=k  }{ \sum}} u_{pl} u_{qm}  c_{lmr}. \notag
\end{gather}  
Note that since we restrict attention to initial conditions for which the unresolved variables are zero and the projection sets the unresolved variables to zero, the quantity $Pe^{tL}QLu_{0kr}$ can be computed through the evolution of the full system \eqref{burgersodemz}.

\subsection{Numerical results}\label{numerical}

In this section we present numerical results for the reduced model of the viscous Burgers equation with viscosity coefficient $\nu = 0.03.$ The solution of the full system was computed with $N=196$ Fourier modes ($F=[-98,97]$) and the first 7 Legendre polynomials ($M=7$). The first 7 Legendre polynomials were enough to obtain converged statistics for the full system. The full system was solved with the modified Euler method with $\delta t = 0.001.$

The reduced model uses $N=196$ Fourier modes but only the first two Legendre polynomials, so $\Lambda=2.$ It was solved using the modified Euler method with $\delta t = 0.001.$ The parameter $t_0$ needed for the evolution of the memory term was found to be 0.3783 through the procedure described in Section \ref{mz_optimal}. 


\begin{figure}
\centering
\epsfig{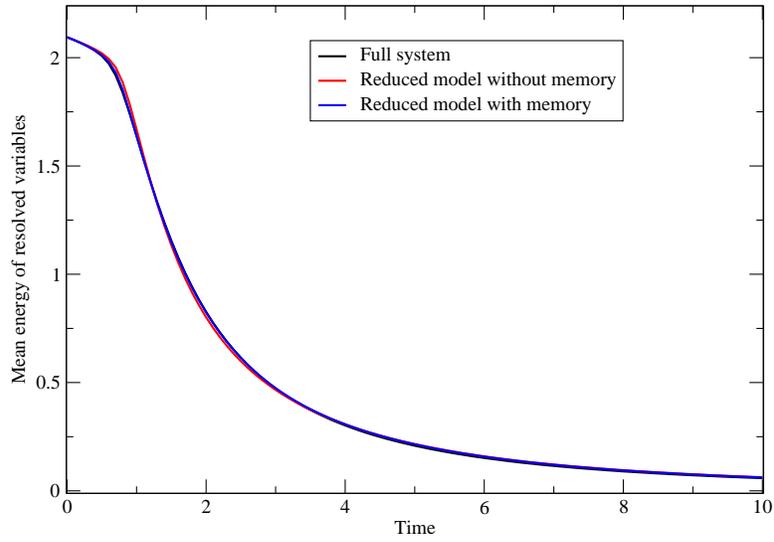}
\caption{Evolution of the mean of the energy of the solution using only the first two Legendre polynomials.}
\label{plot_initial_energy_mean}
\end{figure}

\begin{figure}
\centering
\epsfig{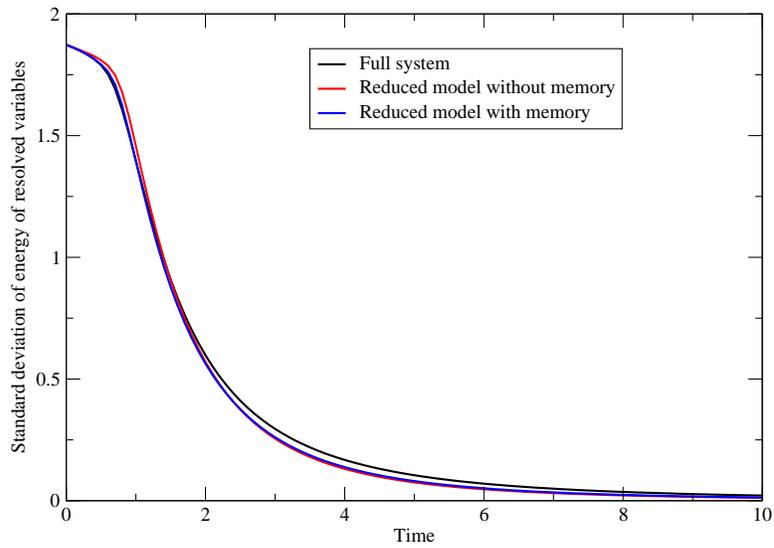}
\caption{Evolution of the standard deviation of the energy of the solution using only the first two Legendre polynomials.}
\label{plot_initial_energy_stdev}
\end{figure}

Figure \ref{plot_initial_energy_mean}  shows the evolution of the mean energy of the solution 
$$\mathbb{E}[E(t)]=\frac{1}{2}  \sum_{k \in F} \sum_{r=0}^1 2\pi |u_{kr}(t)|^2 \frac{1}{2r+1}$$
as computed from the full system (with $M=7$ Legendre polynomials), the MZ reduced model with $\Lambda=2$ {\it without} memory (keeping only the Markovian term) and the MZ reduced model with $\Lambda=2$ {\it with} memory. Figure \ref{plot_initial_energy_stdev} shows the evolution of the standard deviation of the energy of the solution. The variance of the energy is given by
$$Var[E(t)]=\frac{1}{4}  \sum_{k_1, k_2 \in F}   \sum_{r_1,\ldots,r_4=0}^1 (2\pi)^2 u_{k_1r_1} u_{k_1r_2}^*u_{k_2r_3} u_{k_2r_4}^*d_{r_1r_2r_3r_4}-\{ \mathbb{E}[E(t)]\}^2,$$
where
$$d_{r_1r_2r_3r_4}=\int_{-1}^1L_{r_1}(\xi)L_{r_2}(\xi)L_{r_3}(\xi)L_{r_4}(\xi) \frac{1}{2} d\xi.$$
The reduced model performs equally well with or without memory. Of course, the reduced model with memory is slower than the reduced model without memory. However, the reduced model with memory is still about 4 times faster than the the full system.

\begin{figure}
\centering
\epsfig{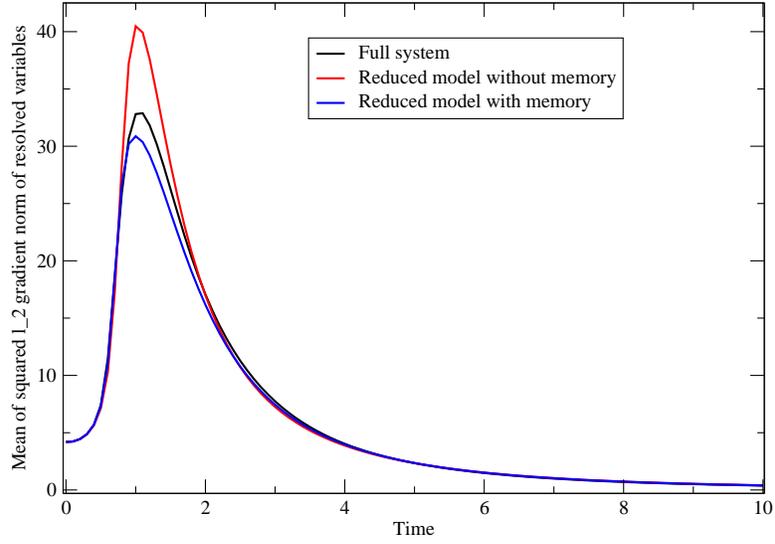}
\caption{Evolution of the mean of the squared $l_2$ norm of the gradient of the solution calculated using only the first two Legendre polynomials.}
\label{plot_initial_gradient_mean}
\end{figure}

\begin{figure}
\centering
\epsfig{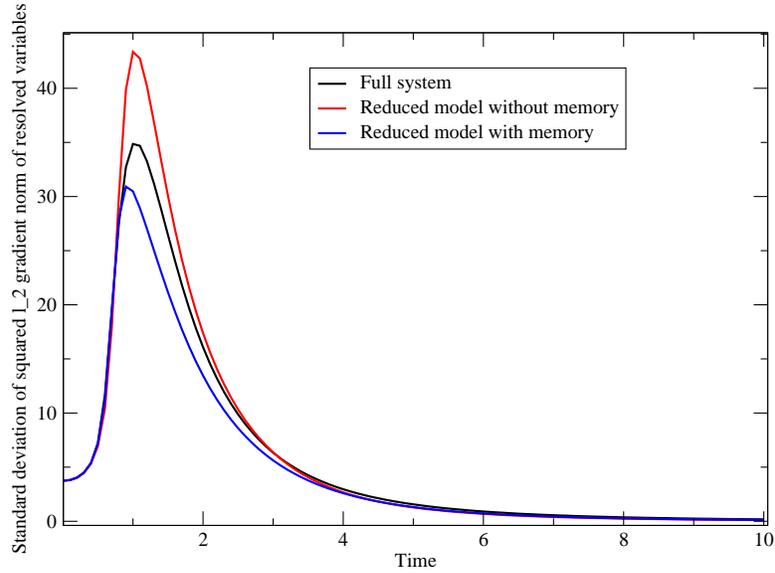}
\caption{Evolution of the standard deviation of the squared $l_2$ norm of the gradient of the solution calculated using only the first two Legendre polynomials.}
\label{plot_initial_gradient_stdev}
\end{figure}

Figure \ref{plot_initial_gradient_mean} shows the evolution of the mean squared $l_2$ norm of the gradient of the solution 
$$\mathbb{E}[G(t)]= \sum_{k \in F} \sum_{r=0}^1 2\pi k^2 |u_{kr}(t)|^2 \frac{1}{2r+1}$$
as computed from the full system (with $M=7$ Legendre polynomials), the MZ reduced model with $\Lambda=1$ {\it without} memory (keeping only the Markovian term) and the MZ reduced model with $\Lambda=1$ {\it with} memory. Figure \ref{plot_initial_gradient_stdev}  shows the evolution of the standard deviation. The variance is given by $$Var[G(t)]= \sum_{k_1, k_2 \in F}   \sum_{r_1,\ldots,r_4=0}^1 (2\pi)^2 k_1^2 k_2^2u_{k_1r_1} u_{k_1r_2}^*u_{k_2r_3} u_{k_2r_4}^*d_{r_1r_2r_3r_4}-\{ \mathbb{E}[G(t)]\}^2.$$
It is obvious from the figures that the inclusion of the memory term improves the performance of the reduced model. Recall that the solution of Burgers equation is a contraction \cite{lax}. Eventually, the complete description of the uncertainty caused by the uncertainty in the initial condition requires only a few polynomial chaos expansion coefficients. This happens at a time scale that is dictated by the magnitude of the viscosity coefficient. That is why for long times the reduced model with and without memory have comparable behavior to that of the full system. However, for short times, the inclusion of the memory term does make a difference because information from the higher chaos expansion coefficients is needed. The higher chaos expansion coefficients will have a more prolonged contribution for systems that possess unstable modes. In such cases, the inclusion of the memory term becomes imperative for short as well long times. Results for such cases will be presented elsewhere.

\section{Discussion and future work}\label{discussion}

We have presented the application of the Mori-Zwanzig formalism to the construction of reduced models for systems of differential equations resulting from polynomial chaos expansions of solutions of differential equations with initial condition uncertainty. In particular, we presented a way that the reduced model can be reformulated so that instead of integro-differential equations one has to solve differential equations. The problem that arises in any reduced model with memory is to compute the length of the memory. We have presented an algorithm which allows the estimation of the memory length on the fly. This means that one starts evolving the full system and use it to estimate the reduced model parameters. Once this is achieved, the simulation continues by evolving only the reduced model with the necessary parameters set equal to their estimated values from the first part of the algorithm. The construction presented here can also be applied to the problem of constructing reduced models for systems forced by random noise \cite{hou}. 

In the current work we have focus on the simplest form of the Markovian reformulation of the Mori-Zwanzig formalism which keeps only the first term in the expansion of the memory. In \cite{s12}, where we dealt with the problem of constructing reduced models when the original system exhibits parametric uncertainty we presented reduced models which used the first two terms in the memory expansion. However, the estimation of the necessary parameters for the reduced model was addressed in a different way. The method presented in \cite{s12} cannot be applied in the current case and this is the reason we had to devise an alternative approach. Results from the application of the current approach to reduced models where one keeps more terms in the memory expansion will be presented in a future publication.    

Finally, we mention that one can construct models which effect reduction both for the variables needed to describe uncertainty and the number of variables needed to describe the system for one realization of the uncertainty sources. This two-level reduction is imperative in situations where solving even for one realization of the uncertainty sources is very expensive e.g. atmospheric flows, fluid structure interactions.

\section*{Acknowledgements}
I would like to thank Prof. G. Karniadakis and Dr. D. Venturi for useful discussions and comments.


\begin{thebibliography}{10} 

\bibitem{barreira} {\sc Barreira L. and Valls C.}, {\em Ordinary Differential Equations: Qualitative Theory}, American Mathematical Society, 2012.

\bibitem{CHK00}
{\sc Chorin, A.J., Hald, O.H. and Kupferman, R.},
{\em Optimal prediction and the Mori-Zwanzig representation of irreversible
processes}, Proc. Nat. Acad. Sci. USA 97 (2000) pp. 2968-2973.

\bibitem{CHK3}
{\sc Chorin, A.J., Hald, O.H. and Kupferman, R.}, {\em Optimal prediction with memory}, Physica D 166 (2002) pp. 239-257.

\bibitem{CS05}
{\sc Chorin, A.J. and Stinis, P.}, {\em Problem reduction, renormalization and memory},
Comm. App. Math. Comp. Sci. 1 (2005) pp. 1-27.

\bibitem{evans}
{\sc Evans L.C.}, {\em Partial Differential Equations}, Second Edition, American Mathematical Society, 2010.

\bibitem{ghanem}
{\sc Ghanem R. and Spanos P.D.},{\em Stochastic finite elements: a spectral approach}, Springer-Verlag, 1998.

\bibitem{givon}
{\sc Givon, D., Kupferman, R. and Stuart, A.}, {\em Extracting macroscopic 
dynamics: model problems and algorithms}, Nonlinearity 17 (2004) pp. R55-R127.

\bibitem{gupta}
{\sc Gupta M. and Narasimhan S.G.}, {\em Legendre polynomials Triple Product Integral and lower-degree approximation of polynomials using Chebyshev polynomials}, Technical Report - CMU-RI-TR-07-22, Carnegie Mellon, 2007.

\bibitem{hou}
{\sc Hou T.Y., Luo W., Rozovskii B. and Zhou H.M.}, {\em Wiener Chaos Expansions and Numerical Solutions of Randomly Forced Equations of Fluid Mechanics}, J. Comput. Phys. 216 (2006) pp. 687-706. 

\bibitem{lax}
{\sc Lax, P.D.}, {\em Hyperbolic Systems of Conservation Laws and the Mathematical Theory of
Shock Waves}, SIAM Publications, Philadelphia, 1972.

\bibitem{leonenko}
{\sc Leonenko G. and Phillips T.}, {\em On the solution of the Fokker-Planck equation using a high-order reduced basis
approximation}, Comput. Methods Appl. Mech. Engrg., 199(1-4) (2009) pp. 158-168.

\bibitem{ma}
{\sc Ma X. and Zabaras N.}, {\em An adaptive hierarchical sparse grid collocation method for the solution of stochastic
differential equations}, J. Comput. Phys., 228 (2009) pp. 3084-3113.

\bibitem{nouy}
{\sc Nouy A. and Le Ma\^itre O. P.}, {\em Generalized spectral decomposition for stochastic nonlinear problems}, J. Comput.
Phys., 228 (2009) pp. 202-235.

\bibitem{s12}
{\sc Stinis P.}, {\em Mori-Zwanzig reduced models for uncertainty quantification I: Parametric uncertainty}, arXiv:1211.4285v1.

\bibitem{venturi}
{\sc Venturi D.}, {\em A fully symmetric nonlinear biorthogonal decomposition theory for random fields}, Physica D,
240(4-5) (2011) pp. 415-425.

\bibitem{wan}
{\sc Wan X. and Karniadakis G. E.}, {\em Multi-element generalized polynomial chaos for arbitrary probability measures},
SIAM J. Sci. Comput., 28(3) (2006) pp. 901-928.

\end{thebibliography}
\end{document}